\documentclass[12pt]{amsart}
\usepackage{amssymb}
\newtheorem{thm}{Theorem}
\newtheorem{lemma}[thm]{Lemma}
\newtheorem{prop}[thm]{Proposition}
\newtheorem{cor}[thm]{Corollary}
\newtheorem{rmk}[thm]{Remark}

\newcommand{\mysection}[1]{\section{#1}\setcounter{equation}{0}}
\def\pf{{\textit{Proof :} }}
\def\qed{\hfill{Q.E.D.}\smallskip}

\newcommand{\ver}{\varepsilon}

\newcommand{\conn}{\nabla}
\newcommand{\pconn}{\nabla ^{\partial M}}
\newcommand{\pD}{D ^{\partial M}}
\newcommand{\D}{\Delta}
\newcommand{\s}{\mathbb{S}}

\newcommand{\p}{\partial}

\let\phi=\varphi
\begin{document}

\title [Dirac Operator on Manifolds with Boundary]{Eigenvalues  of the Dirac
Operator \\on Manifolds with Boundary}

\author{Oussama Hijazi}
\address[Hijazi]{Institut {\'E}lie Cartan\\
Universit{\'e} Henri Poincar{\'e}, Nancy I\\
B.P. 239\\
 54506 Vand\oe uvre-L{\`e}s-Nancy Cedex, France}
\email{hijazi@iecn.u-nancy.fr}

\author{Sebasti{\'a}n Montiel}
\address[Montiel]{Departamento de Geometr{\'\i}a y Topolog{\'\i}a\\
Universidad de Granada\\
18071 Granada \\
Spain}
\email{smontiel@goliat.ugr.es}

\author{Xiao Zhang}
\address[Zhang]{Institute of Mathematics\\
Academy of
Mathematics and Systems Sciences, Chinese Academy of
Sciences\\
Beijing 100080, P.R. China}
\email{xzhang@math08.math.ac.cn}



\begin{abstract} Under standard local boundary conditions or certain global 
APS boundary conditions, we get lower bounds for the eigenvalues of the Dirac 
operator on compact spin manifolds with boundary. Limiting cases are 
characterized by the existence of real Killing spinors and the minimality of 
the boundary.
\end{abstract}

\thanks{Research of S.M. is partially
supported by a DGICYT grant No. PB97-0785. Research of X.Z. is
partially supported by the Chinese NSF and mathematical physics
program of CAS}

\thanks{This work is partially done during the visit of the last two
authors to the Institut {\'E}lie Cartan, Universit{\'e} Henri
Poincar{\'e}, Nancy 1.  They would like to thank the institute
for its hospitality.}


\date{August 17, 2000}

\maketitle \pagenumbering{arabic}

\mysection{Introduction}
It is well known that the spectrum of the Dirac operator on closed spin
manifolds detects subtle information on the geometry and the topology of 
such manifolds (see for example \cite{BFGK,BHMM}).

In \cite{HZ1}, basic properties of the hypersurface Dirac operator
are established. This hypersurface Dirac operator appears as the boundary 
term in
the integral Schr{\"o}dinger-Lichnerowicz formula (\ref{ili}) for compact 
spin manifolds with compact boundary. In fact, the hypersurface Dirac 
operator is, up to a zero order operator,
the intrinsic Dirac operator of the boundary. 

In this paper, we examine the classical local boundary conditions
and certain Atiyah-Patodi-Singer boundary conditions for the Dirac operator.
Here, the spectral resolution of the intrinsic Dirac operator of the boundary 
is used to define the APS boundary conditions. We first prove 
self-adjointness and ellipticity of such conditions.

Then,  systematic use of the modified Levi-Civit{\`a} connections,
introduced in \cite{Fr,Hi3,Zh1,FK,HZ1}, is made. Under appropriate curvature
assumptions, these modified connections  combined
with formula (\ref{ili}), yield to the corresponding estimates for 
compact spin manifolds with boundary. The limiting cases are then studied.

Such estimates are obtained in Sections \ref{lb} and \ref{clb}. In Section 
\ref{lb} we consider both the local and the above mentioned APS boundary 
conditions. We first introduce the modified connection (\ref{gcm})
which allows to establish a Friedrich's type inequality,
in case the mean curvature of the boundary is nonnegative. Under the local
boundary conditions, the limiting
case is then charaterized by the existence of a Killing spinor on the 
compact manifold with minimal boundary (see (\ref{est1})). Then the 
energy-momentum tensor is used to define the modified connection (\ref{emtm}),
from which one can deduce inequality (\ref{est2}).

Finally, in section \ref{clb}, under the 
local boundary conditions, the conformal aspect is examined. For example, 
generalizations of the conformal lower bounds 
in \cite{Hi1,Hi3} are obtained (see Remark \ref{Hi1Hi3}). 

\mysection{The elliptic boundary conditions}
Let $M$ be an $n$-dimensional Riemannian spin manifold with boundary $\p M$
endowed with its induced Riemannian and spin structures.
Denote by $\s$ the spinor bundle of $M$. Let $\conn$ (resp. $\conn ^{\p M}$)
be the Levi-Civit{\`a} connection of $M$ (resp. $\p M$) and denote by the same
symbol their
corresponding lift to the spinor bundle $\s$. Consider the Dirac
operator $D$ of $M$ defined by $\conn $ on $\s$. It is known
\cite{LM} that there exists a positive definite Hermitian metric
on $\s $ which satisfies, for any covector field $X^*\in \Gamma
(T^*M)$, and any spinor fields $\phi , \psi \in \Gamma (\s) $, the
relation
\begin{eqnarray}
(X^*\cdot\phi,\; X^*\cdot\psi)=|X^*|^2(\phi,\;\psi),\label{ip}
\end{eqnarray}
where ``$\cdot $'' denotes  Clifford multiplication. The connection
$\conn $ is compatible with the metric $(\; , \;)$. Fix a point
$p\in \p M$ and an orthonormal basis $\{ e _{\alpha}\}$ of $T_p M$
with $e_0$ the outward normal to $\p M$ and $e_i$ tangent to $\p M$
such that for $1\leq i,j \leq n$ $$(\pconn  _i e_j)_p = (\conn _0
e_j)_p =0.$$ Let $\{e^{\alpha}\}$ be the dual coframe. Then, for
$1\leq i,j \leq n$
\begin{eqnarray*}
(\conn _i e^j)_p &=& -h_{ij}\;e^0,\\ (\conn _i e^0)_p &=&
h_{ij}\;e^j,
\end{eqnarray*}
where $h_{ij}$ = $(\conn_i e_0,e_j)$ are the components of the
second fundamental form at $p$, and we have
\begin{eqnarray}
\conn _{i} = \pconn _{i}+\frac{1}{2} \;h_{ij}\; e^0\cdot
e^j\cdot\;.  \label{rtc2}
\end{eqnarray}
Let $H =\sum h _{ii}$ be the mean curvature of $M$. In the
above notation, the standard sphere $S ^n _r =\p B ^{n+1} _r $
has positive mean curvature $H = \frac{n}{r}$. By (\ref{ip}),
$(e^0\cdot e^j\cdot \phi,\psi)=(\phi, e^j\cdot e^0\cdot \psi)$.
Therefore (\ref {rtc2}) implies
\begin{eqnarray*}
d(\phi,\;\psi)*e_i & = &\Big((\conn _{i} \phi,\psi)+(\phi, \conn
_{i} \psi)\Big)*1\\ & = &\Big((\pconn _{i}\phi,\psi)+(\phi,\pconn
_{i} \psi)\Big)*1.
\end{eqnarray*}
Hence the connection $\pconn $ is also compatible with the metric
$(\; , \;)$. Denote by $\pD $ the Dirac operator of $\p M$. In the
above orthonormal coframe $\{e^i\}$ of $M$, $\pD = e^i\cdot\pconn
_{i}$. Thus $\pD$ is self-adjoint with respect to the metric $(\;
, \;)$. The relation (\ref{rtc2}) implies that $$\pconn _i (e ^0\cdot  \phi) =
e^0 \cdot \pconn _i \phi.$$ Hence $$\pD (e ^0\cdot  \phi ) = -e
^0\cdot \pD \phi.$$

Consider the integral form of the Schr{\"o}dinger--Lichnerowicz formula for
a compact manifold with compact boundary
\begin{eqnarray}
\int _{\p M} (\phi,\;e ^0 \cdot \pD \phi)& - &\frac{1}{2} \int _{\p
M} H |\phi |^2 = \nonumber\\ & &\int _M |\conn \phi |^2
+\frac{R}{4}|\phi |^2 -|D
\phi |^2 .   \label{ili}
\end{eqnarray}

It is well-known that there are basically two types of elliptic
boundary conditions for the Dirac operator: The local boundary
condition and the (global) Atiyah-Patodi-Singer (APS) boundary
condition. Such boundary conditions  are used in the positive mass
theorem for black holes, Penrose conjecture in general relativity
and the index theory in topology \cite{GHHP, G, He1, He2, Zh2}.
The APS boundary condition exists on any spin manifold with
boundary, while the local boundary condition requires certain
additional structures on manifolds such as the existence of a
Lorentzian structure or a chirality operator, etc \cite{FS, GHHP,
He2}. Now we shall show that the local boundary condition exists
on certain spin manifolds with a ``boundary chirality operator''.

An operator $\Gamma $ defined on $C ^\infty (\p M, \s | _{\p M} )$
is said to be a boundary chirality operator if it satisfies the 
following conditions:
\begin{eqnarray}
\Gamma ^2 & = & Id, \\ \conn ^{\partial M} _ {e _i } \Gamma &=&0,\\ e ^ 0
\cdot \Gamma  &=& - \Gamma \cdot e ^0 ,\\e ^i \cdot
\Gamma  &=& \Gamma \cdot e ^i ,\\ (\Gamma \cdot \phi,
\Gamma \cdot \psi) &=& (\phi, \psi).
\end{eqnarray}

If $M$ is a spacelike hypersurface of a spacetime manifold with timelike
covector $T$, then we can let $\Gamma = T \cdot e ^0$, where $e
^0$ is the normal covector on $\partial M$. 

Recall that (see \cite{FS} for example), an operator $F$ defined on
$C ^\infty (M, \s )$ is {\it called} a chirality operator on $M$ if for 
all $X^*\in\Gamma (T^*M)$, and any spinor fields $\phi , \psi \in \Gamma (\s)$,
one has
\begin{eqnarray*}
F ^2 & = & Id, \\ 
\conn _ {X} F &=&0,\\ 
X^*\cdot F  &=& - F \cdot X^*,\\ 
(F\cdot \phi, F\cdot \psi) &=& (\phi, \psi).
\end{eqnarray*}
Note that, such an operator exists if the spin manifold $M$ is even 
dimensional. It is easy to see that if $M$ has a chirality operator $F$, then
$\Gamma = F | _{\partial M}\cdot e ^0$ is a boundary chirality operator. But, 
in general, it is not
known whether the local boundary condition exists. In this paper, 
we consider 
the following boundary conditions:\\

\begin{itemize}
\item {\bf The local boundary condition}: \\
As the eigenvalues of the chirality operator $\Gamma$ 
are $\pm 1$, the corresponding eigenspaces
\begin{eqnarray*}
\Gamma ^{loc} _{+} &=&\Big\{ \phi \in C ^\infty (\p M, \s | _{\p
M} ),\;\; \Gamma \cdot \phi = \phi \Big\},\\ \Gamma ^{loc} _{-}
&=&\Big\{ \phi \in C ^\infty (\p M, \s | _{\p M} ),\;\;\Gamma
\cdot\phi = -\phi \Big\}
\end{eqnarray*}
provide local boundary conditions.\\

\item {\bf The APS type boundary condition}: \\
The operator $e ^0
\cdot \pD $ is self-adjoint with respect to the induced metric
$(\;\;,\;\;)$ on $\p M$. Therefore it has a discrete (real)
spectrum. Let $(\phi _k) _{k \in \mathbb{N} }$ be the spectral
resolution of $e ^0 \cdot \pD$, i.e., $e ^0 \cdot \pD \phi _k =
\lambda _k \phi _k $, and consider the corresponding 
$L^2$-orthogonal subspaces $\Gamma ^{APS} _{\pm}$ spanned by the 
positive and negative eigenspaces of $e ^0 \cdot\pD $, i.e.,
\begin{eqnarray*}
\Gamma ^{APS} _{+} &=&\Big\{ \phi \in C ^\infty (\p M, \s | _{\p
M} ),\;\; \phi = \sum _{\lambda _k > 0} c _k \phi _k\Big\},\\
\Gamma ^{APS} _{-} &=&\Big\{ \phi \in C ^\infty (\p M, \s | _{\p
M} ), \;\;\phi = \sum _{\lambda _k < 0} c _k \phi _k\Big\}.
\end{eqnarray*}
We will
consider the APS type boundary conditions corresponding to the 
projections onto these subspaces.
Recall that the original
Atiyah-Patodi-Singer (APS) boundary condition refers to the 
spectral resolution of $e^0\cdot \pD-H/2 $ instead of  $e^0\cdot\pD$ 
(see for example \cite{FS}).
\end{itemize}

Note that if $\phi,\psi \in \Gamma _{\pm} ^{loc}$, then
\begin{eqnarray*}
(e ^0 \cdot \phi, \psi) &=& (e ^0 \cdot \Gamma \cdot \phi, \Gamma
\cdot \psi)\\&=& -(\Gamma \cdot e ^0 \cdot \phi, \Gamma \cdot
\psi)\\&=&-(e ^0 \cdot \phi, \psi),
\end{eqnarray*}
hence $(e ^0 \cdot \phi, \psi) =0$. On the other hand, if $\phi,\psi
\in \Gamma _{\pm} ^{APS}$, then \break $e ^0 \cdot \phi \in \Gamma _{\mp} ^
{APS}$, therefore
\begin{eqnarray*}
\int _{\p M} (e ^0\cdot \phi,\psi)=0.
\end{eqnarray*}
These facts imply that the Dirac operator $D$ is self-adjoint
under either the local boundary conditions or the APS boundary
negative and positive conditions. Morevoer, it has real eigenvalues. Now 
we define the $H ^k$
Sobolev norm by $$
\parallel \phi
\parallel ^2 _{H ^k} = \sum _{|\alpha |=k }
\parallel \conn ^ \alpha \phi \parallel ^2 _ {L ^2} +
\parallel \phi \parallel ^2 _{H ^{k-1}},$$
where $ \alpha $ is a multi-index.
\begin{prop}
Under the local boundary condition $\phi \in \Gamma ^{loc} _{\pm}$
or the APS boundary condition $\phi \in \Gamma ^{APS} _{-}$, the
Dirac operator $D$ satisfies elliptic estimates: For any
$k \geq 1$, $\delta >0$, there exists $C _{k,\delta} $ such that
\begin{eqnarray}
\parallel \phi \parallel ^2 _ {H ^k} \leq (1+\delta )
\parallel D \phi \parallel ^2 _{L ^2} + C _{k, \delta }
\parallel \phi \parallel ^2 _ {H ^{k-1}}.  \label{heloc}
\end{eqnarray}
\end{prop}
\pf Note that for any $\phi \in \Gamma ^{loc} _{+}$ or $\phi \in
\Gamma ^{loc} _{-}$, $D ^{\partial M} (\Gamma \cdot \phi) = \Gamma
\cdot D ^{\partial M} \phi$, thus
\begin{eqnarray*}
(\phi, e ^0 \cdot \pD \phi) &=&\Big(\Gamma \cdot \phi, e ^0 \cdot
\pD (\Gamma \cdot \phi)\Big)\\&=& (\Gamma \cdot \phi, e ^0 \cdot
\Gamma \cdot \pD \phi)\\&=&-(\phi, e ^0 \cdot \pD \phi).
\end{eqnarray*}
Therefore $(\phi, e ^0 \cdot \pD \phi)=0$. If $\phi \in \Gamma
^{APS} _{-}$, then $$\int _{\p M} (\phi, e ^0 \cdot \pD \phi)=\sum
_{\lambda _k < 0} |c _k |^2 \lambda _k \le 0.$$ By the
Ehrling-Gagliardo-Nirenberg inequality \cite{A}, for each $\ver >0$,
there exists a constant  $C _\ver >0$ such that
\begin{eqnarray*}
\parallel \phi \parallel ^2 _ {L ^2 ({\p M})} \leq \ver
\parallel \phi \parallel ^2 _{H ^1} + C _\ver
\parallel \phi \parallel ^2 _{L ^2},
\end{eqnarray*}
thus (\ref{ili}) implies
\begin{eqnarray}
\parallel \phi \parallel ^2 _ {H ^1} \leq (1+\delta )
\parallel D \phi \parallel ^2 _{L ^2} + C _\delta
\parallel \phi \parallel ^2 _{L ^2}.  \label{eloc}
\end{eqnarray}
Then a standard argument gives (\ref {heloc}).\qed

The following corollary is a direct consequence of the Sobolev
embedding theorem $$\parallel \phi \parallel ^2 _ {C ^k} \leq
C \parallel \phi \parallel ^2 _ {H ^{k +\frac{n}{2}}}$$
\begin{cor}
Any eigenspinor of the Dirac operator which satisfies either the local
boundary condition $\phi \in \Gamma ^{loc} _{\pm}$ or the
(negative) APS boundary condition $\phi \in \Gamma ^{APS} _{-}$ is
smooth.
\end{cor}
\mysection {Lower Bounds for the Eigenvalues}\label{lb}
In this section, we adapt the arguments used in \cite{HZ1} to the case of spin
compact manifolds with boundary. In particular, we get 
generalizations of basic inequalities on the eigenvalues of the Dirac operator
$D$ under the local 
boundary conditions $\Gamma ^{loc} _\pm $ or the negative APS boundary 
condition $\Gamma_-^{APS}$. 
For this, we use the integral identity (\ref{ili}) together with
an appropriate modification of the Levi-Civit{\`a} connection. 

Let $D \phi = \lambda \phi$, where $\lambda $ is a real constant or
a real function. For any real functions $a$ and $u$, we define
\begin{eqnarray}\label{gcm}
\conn ^{a,u} _i =\conn _i +a \conn _i u +\frac{a}{n} \conn _j u \;
e ^i \cdot e ^j \cdot +\frac{\lambda }{n} e^i \cdot .
\end{eqnarray}
Then
\begin{eqnarray*}
| \conn ^{a,u} \phi | ^2 &=&| \conn \phi | ^2+\frac{\lambda ^2}{n}
|\phi | ^2 + a ^2 (1-\frac{1}{n}) |du |^2 |\phi |^2\\&& +2 a \conn
_i u \;\Re (\conn _i \phi, \; \phi)+\frac{2\lambda }{n} \;\Re (\conn
_i \phi, \;e ^i \cdot \phi)\\ &=&| \conn \phi | ^2 -\frac{\lambda
^2}{n} |\phi |^2+ a ^2 (1-\frac{1}{n}) |du |^2 |\phi |^2 +a \conn
_i u \conn _i |\phi |^2.
\end{eqnarray*}
Define the functions $R_{a,u}$ by
\begin{eqnarray}\label{Rau}
R_{a,u} =R -4a \D u +4 \conn a \conn u -4 (1-\frac{1}{n})a ^2 |du
|^2
\end{eqnarray}
where $\D $ is the positive scalar Laplacian. Then we have
\begin{eqnarray*}
\int _M | \conn ^{a,u} \phi | ^2 &=&\int _M | \conn \phi | ^2
-\frac{\lambda ^2}{n} |\phi |^2 -\Big(\frac{R _{a,u}}{4}
-\frac{R}{4}\Big) |\phi |^2\\
&& +\int _{\p M} a\;du(e _0) |\phi |^2.
\end{eqnarray*}
Therefore (\ref{ili}) yields
\begin{eqnarray}
\int _M |\conn ^{a,u} \phi| ^2 &=& \int _M
\Big[(1-\frac{1}{n})\lambda ^2 -\frac{R _{a,u}}{4} \Big] |\phi |^2
 \nonumber\\ &&+\int _{\p M}
(\phi,\;e ^0 \cdot \pD \phi) +\Big[a\;du(e _0)-\frac{H}{2}\Big]
|\phi |^2 . \label{eq1}
\end{eqnarray}

Now we generalize Lemma 2.3 in \cite{FK} to the case where $a$ is a
real function.
\begin{lemma}\label{lemma1}
Suppose there exist a spinor field $\phi \in \Gamma (\s )$, a real 
number $\lambda $ and  a real
functions $a$ and $u$ on $M$ such that for all $i$, $1\le i\le n$,
\begin{eqnarray}
\conn _i \phi = - \frac{\lambda }{n} e ^i \cdot \phi - a \conn _i u
\phi -\frac{a}{n} \conn _j u e ^i \cdot e^j \cdot \phi . \label{ieq1}
\end{eqnarray}
Then $\phi$ is a real Killing spinor, i.e., either $a=0$ or $du=0$. In particular, the manifold is Einstein.
\end{lemma}
\pf First, observe that (\ref{ieq1}) implies  $D \phi = \lambda \phi$. By
the Ricci  identity (see \cite{FK}), we have
\begin{eqnarray*}
\frac{1}{2} R _{ij} e ^i \cdot e ^j \cdot \phi &=&e ^i \cdot D (\conn
_i \phi) - D ^2 \phi \\
&=&e ^i \cdot e ^j \cdot \conn _j \big(
-\frac{\lambda }{n} e ^i \cdot - a \conn _i u
 -\frac{a}{n} \conn _k u e ^i\cdot e^k\cdot\big)\phi - \lambda
^2 \phi \\
&=&\frac{\lambda }{n} e ^i \cdot ( e ^i \cdot e ^j \cdot + 2 \delta
_{ij}) \conn _j \phi \\
& & - \lambda ^2 \phi -du \cdot da \cdot \phi +a
\triangle u \phi -a \lambda du \cdot \phi\\
& &+\frac{1}{n} e ^i \cdot (e ^i \cdot e ^j \cdot + 2\delta _{ij} ) e
^k \cdot \Big(\conn _j a \conn _k u \phi \\
& & + a \conn _j \conn _k u \phi
+ a \conn _k u \conn _j \phi \Big)\\
&=&\Big(\frac{2(1-n)}{n} \lambda ^2  +\frac{2a}{n} \triangle u
-\frac{2(2-n)}{n}\conn a \conn u \Big)\phi \\
&& -\frac{2}{n} du \cdot da
\cdot \phi +\frac{4 a \lambda}{n ^2} du \cdot \phi.
\end{eqnarray*}
This implies either $a=0$ or $du=0$.\qed

By (\ref {eq1}) and Lemma \ref {lemma1}, we obtain
\begin{thm}
Let $M^n $ be a compact Riemannian spin manifold of dimension $n \ge 2$, with
boundary $\p M$, and let $\lambda $ be any eigenvalue of $D$ under
either the local boundary condition $\Gamma ^{loc} _\pm $ or the
(negative) APS boundary condition $\Gamma ^{APS} _{-}$. If there
exist real functions $a$, $u$ on $M$ such that $$H \geq 2a \;du(e
_0) $$ on $\p M$, where $H$ is the mean curvature of $\p M$, then
\begin{equation}
\lambda ^2 \geq \frac{n}{4(n-1)}\sup _{a , u} \inf _M R
_{a,u},\label{est1}
\end{equation}
where $R_{a,u}$ is given in (\ref{Rau}).
In the limiting case with the local boundary conditions,  the associated eigenspinor is a real Killing spinor and $\p M$ is minimal. 
\end{thm}

Note that by \cite{HMZ}, under the APS
boundary conditions equality in (\ref{est1}) could not hold.

Now we make use of  the energy-momentum tensor (see \cite{Hi3}) to get
lower bounds for the eigenvalues of $D$. For any spinor field $\phi$, we
define the associated energy momentum 2-tensor $Q _{\phi }$ on the
complement of its zero set by,
\begin{eqnarray}
Q_{\phi, ij}=\frac{1}{2}\;\Re (e ^i\cdot\conn _j \phi + e ^j
\cdot\conn _i \phi\;, \phi/|\phi |^2). \label{gq}
\end{eqnarray}
If $\phi $ is an eigenspinor of $D$, the tensor $Q _{\phi}$ is
well-defined in the sense of distribution. Let
\begin{eqnarray}\label{emtm}
\conn ^{Q,a,u} _i =\conn _i +a \conn _i u +\frac{a}{n} \conn _j u
\; e ^i \cdot e ^j \cdot +Q _{\phi, ij} e ^j \cdot.
\end{eqnarray}
It is easy to prove that (see \cite{HZ1})
\begin{eqnarray*}
| \conn ^{Q,a,u} \phi | ^2 =| \conn \phi | ^2 -|Q _\phi |^2 |\phi |^2+ a
^2 (1-\frac{1}{n}) |du |^2 |\phi |^2 +a \conn _i u \conn _i |\phi
|^2.
\end{eqnarray*}
Therefore
\begin{eqnarray}
\int _M |\conn ^{Q,a,u} \phi| ^2 &=& \int _M \Big[\lambda ^2
-\Big(\frac{R _{a,u}}{4} +|Q _\phi |^2\Big)\Big] |\phi |^2
\nonumber\\ &&+\int _{\p M} (\phi,\;e ^0 \cdot \pD
\phi)+\Big[a\;du(e _0)-\frac{H}{2}\Big] |\phi |^2. \label{eq2}
\end{eqnarray}
Thus we have
\begin{thm}
Let $M^n $ be a compact Riemannian spin manifold of dimension $n \ge 2$, with
boundary $\p M$, and let $\lambda $ be any eigenvalue of $D$ under
either the local boundary condition $\Gamma ^{loc} _\pm$ or the
(negative) APS boundary condition $\Gamma ^{APS} _{-}$. If there
exist real functions $a$, $u$ on $M$ such that $$H \geq 2a \;du(e
_0) $$ on $\p M$, where $H$ is the mean curvature of $\p M$, then
\begin{equation}
\lambda ^2 \geq \sup _{a , u}\inf _M \Big(\frac{R _{a,u}}{4} +|Q
_\phi |^2\Big).\label{est2}
\end{equation}
In the limiting case,  one has $H = 2a du(e _0)$ on $\p M$.
\end{thm}
\begin{rmk} Under either the local boundary condition $\Gamma ^{loc} _\pm $ or
the (negative) APS boundary condition $\Gamma ^{APS} _{-} $, assume that
$H\ge 0$. Take
$a =0$ or $u$ constant in (\ref{est1}) and (\ref{est2}),  then one gets
Friedrich's inequality~\cite{Fr}
\begin{eqnarray}
\lambda ^2 \geq \frac{n}{4(n-1)}\inf _M R
\end{eqnarray}
and the following inequality \cite{Hi3}
\begin{eqnarray}
\lambda ^2 \geq \inf _M \Big(\frac{R}{4} +|Q _\phi|^2 \Big).
\end{eqnarray}
\end{rmk}

\mysection {Conformal Lower Bounds}\label{clb}
As in the previous section and under the local 
boundary conditions $\Gamma ^{loc} _\pm $, we show that the conformal 
arguments used in \cite{HZ1} combined with the integral formula 
(\ref{ili})  yield to 
generalizations of all known lower bounds for the eigenvalues of 
the Dirac operator. 

Let $g$ be the metric of $M$. For any real function $u$ 
on $M$,
consider a conformal metric $\bar g = e^{2u} g$. Denote by
$\overline D$ the Dirac operator with respect to this conformal
metric. If $D \phi = \lambda \phi$, then $\overline D\;\overline
\psi = \lambda\; e ^{-u}\; \overline \psi$ where $\overline \psi =
e ^{-\frac{n-1}{2} u}\; \overline \phi$. Note that
\begin{eqnarray*}
\overline{\conn ^{a,u} _{\overline {e _i}}} &=&\overline \conn _
{\overline {e _i}} + a\; e ^{-u}\; \conn _i u +\frac{a}{n} \; e
^{-u}\; \conn _j u \;\overline {e ^i}\; \overline {\cdot}\;
\overline {e ^j}+\frac{\lambda }{n} \; e ^{-u}\; \overline {e
^i}\; \overline {\cdot},\\ \overline \D u &=& - \sum _i e^{-u}
(\conn _{e _i} (e ^{-u} \conn _{e _i} u)) =e ^{-2u} (\D u + |du|
^2),\\ \overline R\; e ^{2u} &=& R + 2(n-1) \D u -(n-1)(n-2)|du|
^2,
\end{eqnarray*}
also, on $\p M$,
\begin{eqnarray*}
 \overline {D ^{\p M}} \Big(\; e ^{-\frac{(n-2)}{2} u}
\;\overline{\phi}\; \Big)&=&e ^{-\frac{n}{2} u} \;\overline {D
^{\p M} \phi},\\ \overline  H &=& e ^{-u} \,\Big( H+ (n-1)\; du (e
_0) \Big).
\end{eqnarray*}

Define the function $\widehat{R}_{a,u}$ by
\begin{eqnarray}\label{HRau}
\widehat{R}_{a,u} &=&R+4(\frac{n-1}{2} -a) \D u +4 \conn a \conn u
\nonumber\\&& - \Big ((n-1)(n-2) +4(2-n)a +4(1-\frac{1}{n})a ^2 \Big ) |du
|^2
\end{eqnarray}
where $\D $ is the positive scalar Laplacian. Then apply
(\ref{eq1}) to the conformal metric $\overline g$, to get
\begin{eqnarray}
\int _M |\; \overline{\widehat{\conn}^{a,u}} \; \overline \psi\;|
^2 _{\bar g}\; \bar v _g &=&\int _M e ^{-u} \;\Big[(1-\frac{1}{n})
\lambda ^2 -\frac{\widehat{R}_{a,u}}{4} \Big]|\phi |^2 v _g
\nonumber\\ &&\!\!\!\!\!\!\!\!+\int _{\p M} \Bigg\{(\overline{\psi},
\;\overline {e
^0}\; \overline{\cdot} \;\overline{\pD }\;\overline {\psi})
_{\overline g}+ \Big[a \;du(\overline {e _0})-\frac{\overline
{H}}{2}\Big] |\overline {\psi} |^2 _{\overline {g}} \Bigg\} v
_{\overline g} \nonumber\\ &=&\int _M e ^{-u}
\;\Big[(1-\frac{1}{n}) \lambda ^2 -\frac{\widehat{R}_{a,u}}{4}
\Big]|\phi |^2 v _g \nonumber\\ &&+\int _{\p M} e ^{-u}
\Bigg\{(\phi,\;e ^0 \cdot \pD
\phi)\nonumber\\&&+\Big[(a-\frac{n-1}{2})\;du(e
_0)-\frac{H}{2}\Big] |\phi |^2 \Bigg\}v _g . \label{eq3}
\end{eqnarray}

Note that $\overline{\widehat{\conn}^{a,u}} \; \overline \psi =0 $
implies
\begin{eqnarray*}
\conn _i \phi = - \frac{\lambda }{n} e ^i \cdot \phi -
(a -\frac{n}{2}) \conn _i u \phi -\frac{1}{n}(a-\frac{n}{2})
 \conn _j u e ^i \cdot e^j \cdot \phi
\end{eqnarray*}
(see \cite{FK}), we thus have either $a =\frac{n}{2}$ or $du =0$ by
Lemma \ref{lemma1}. Thus we obtain:
\begin{thm}
Let $M^n $ be a compact Riemannian spin manifold of dimension $n \ge 2$, with
boundary $\p M$, and let $\lambda $ be any eigenvalue of $D$ under
the local boundary condition $\Gamma ^{loc} _\pm$. If there
exist real functions $a$, $u$ on $M$ such that $$H \geq (2a-n+1)
\;du(e _0) $$ on $\p M$, where $H$ is the mean curvature of $\p
M$, then
\begin{equation}
\lambda ^2 \geq \frac{n}{4(n-1)}\sup _{a , u} \inf _M \widehat{R}
_{a,u},\label{est3}
\end{equation}
where the function $\widehat{R}_{a,u}$ is given in (\ref{HRau}).
In the limiting case,  the associated eigenspinor is a real Killing spinor
and either $H =du (e _0)$ or
$H=0$ on $\p M$.
\end{thm}

Since $\overline Q _{\overline \phi, \bar i\, \bar j} = e ^{-u}\;
{Q _{\phi, ij}}$ under the conformal transformation
$\overline g = e ^{2u} g$, we apply (\ref{eq2}) to the conformal
metric $\overline g$, to get
\begin{eqnarray}
\int _M |\; \overline{\widehat{\conn}^{Q,a,u}} \; \overline
\psi\;| ^2 _{\bar g}\; \bar v _g &=&\int _M e ^{-u} \;\Big[\lambda
^2 -\Big(\frac{\widehat{R}_{a,u}}{4} +|Q _\phi| ^2 \Big)\Big]|\phi
|^2 v _g \nonumber\\ &&+\int _{\p M} e ^{-u} \Bigg\{(\phi,\;e ^0
\cdot \pD \phi)\nonumber\\&&+\Big[(a-\frac{n-1}{2})\;du(e
_0)-\frac{H}{2}\Big] |\phi |^2 \Bigg\}v _g . \label{eq4}
\end{eqnarray}
Thus we have
\begin{thm}
Let $M^n $ be a compact Riemannian spin manifold of dimension $n \ge 2$, with
boundary $\p M$, and let $\lambda $ be any eigenvalue of $D$ under
the local boundary condition $\Gamma ^{loc} _\pm$. If there
exist real functions $a$, $u$ on $M$ such that $$H \geq (2a-n+1)
\;du(e _0) $$ on $\p M$, where $H$ is the mean curvature of $\p
M$, then
\begin{equation}
\lambda ^2 \geq \sup _{a , u}\inf _M \Big(\frac{\widehat{R}
_{a,u}}{4} +|Q _\phi |^2\Big).\label{est4}
\end{equation}
In the limiting case one has $H = (2a-n+1) \;du(e _0)$ on $\p M$.
\end{thm}
\begin{rmk}\label{Hi1Hi3}
If $n \geq 3$,  take $a =0$ and $u=-\frac{2}{n-2}\log h $
in (\ref{est3}) and (\ref{est4}),
where $h$ is a positive eigenfunction of the first eigenvalue
$\mu _1$ of the conformal Laplacian 
$$L := 4 \frac{n-1}{n-2} \triangle +R$$
under the boundary condition
$$dh (e _0)
-\frac{(n-2)H}{2(n-1)}\; h =0.
$$  
Then, one gets the  lower bounds \cite{Hi1, Hi3}
\begin{eqnarray}\label{Mu1}
\lambda ^2 \geq \frac{n}{4(n-1)}\mu _1,
\end{eqnarray}
and
\begin{eqnarray}
\lambda ^2 \geq \inf _M \Big(\frac{\mu _1}{4} +|Q _\phi|^2 \Big)
\end{eqnarray}
under the local boundary condition $\Gamma ^{loc} _{\pm} $. In the
limiting case of (\ref{Mu1}),  the associated eigenspinor is a real Killing spinor and $\p M$ is minimal.
\end{rmk}

\end{document}